\documentclass[a4paper]{article}\nonstopmode
\input amssym.def
\input amssym
\font\teneusm=eusm10
\font\seveneusm=eusm7
\font\fiveeusm=eusm5
\newfam\eusmfam
\textfont\eusmfam=\teneusm
\scriptfont\eusmfam=\seveneusm
\scriptscriptfont\eusmfam=\fiveeusm
\def\script#1{{\fam\eusmfam\relax#1}}

\newtheorem{thm}{Theorem}[section]
\newtheorem{lem}[thm]{Lemma}

\newtheorem{prp}[thm]{Proposition}

\newtheorem{dfn}[thm]{Definition}

\newtheorem{construction}[thm]{Random construction}
\def\dem{{\it Proof.\/ }} 
\def\rem{\refstepcounter{thm}
{\bf Remark \thethm\hskip1ex}}

\def\quest{
\refstepcounter{thm}{\bf Question \thethm}\hskip1ex}

\def\vskipa{\vskip\abovedisplayskip}
\def\vskipb{\vskip\belowdisplayskip}
\def\ie{{\it i.\,e.}~} 
\def\vs{{\it vs.}~} 
\def\af{{\it a fortiori}\/\ }

\def\toll#1{\vtop{\ialign {##\crcr \rightarrowfill \crcr 
\noalign{\kern -1pt \nointerlineskip \vskip2pt} 
$\hfil \scriptstyle{\ #1\ } \hfil $\crcr }}}

\def\({{\fam0\rm (}}
\def\){\/{\fam0\rm )}} 
\def\]{\mathopen]}
\def\[{\mathclose[}
\mathchardef\hook="312C

\def\imp{\Rightarrow}
 
\def\ssi{\Leftrightarrow}

\def\bull{\leavevmode\kern .1ex\vrule height 1ex width .9ex depth
-.1ex \kern .8ex}
\def\ecks{\bull}

\def\eck{\nolinebreak\hspace{\fill}\ecks}
\def\barre{\rule[-2.5pt]{0.6pt}{10pt}}
\def\lmes{\kern0.7pt\barre\kern0.7pt}
\def\rmes{\kern0.7pt\/\barre\kern0.7pt}

\def\bi#1#2{\bigg(\kern-0.5pt{{#1}\atop{#2}}\kern-0.5pt\bigg)}
\def\bip#1#2{\left(\kern-3pt {{#1}\atop{#2}}\kern-2pt \right)}
\def\se{\subseteq}
\def\es{\supseteq}
\def\HM#1{\setbox0=\hbox{#1}\dimen0=\wd0 #1 \kern-\dimen0
\setbox0\hbox{\raise3pt\hbox{$^\frown$}}\advance\dimen0 by -\wd0
\box0\kern\dimen0}

\def\doublesum{\mathop{\sum\sum}\limits}

\def\trace{\mathop{\fam0\rm tr}\nolimits}

\def\P#1{\mathop{\Bbb P}\nolimits
\left[\kern0.5pt{#1}\kern0.5pt\right]}
 
\def\e{\mathop{\fam0 e}\nolimits}
\def\Id{\hbox{\fam0\rm Id}}

\def\eps{\varepsilon}
\def\Row{R}\def\Col{C}

 \def\bip#1#2{\left(\kern-3pt
{{#1}\atop{#2}}\kern-2pt \right)}

\def\DI/{$I\subseteq\N\times\N$}

\def\Ref#1{\hbox{$(\ref{#1})$}}

\def\pstar{^{\scriptscriptstyle(\kern-1pt\lower0.5pt\hbox{$\scriptstyle *$}\kern-1pt)}}
\def\ppstar{^{\scriptscriptstyle(\kern-1pt*\kern-1pt)}}

\def\T{\Bbb T}
\def\N{\Bbb N}

\def\C{\Bbb C}
\def\block#1#2{{\setbox0=\hbox{#1\kern1ex}
           \leftskip=\wd0\parindent=-\wd0\par\leavevmode\box0 #2
           \par}}
\title{Lacunary matrices}
\author{Asma Harcharras, Stefan Neuwirth and Krzysztof Oleszkiewicz}
\date{}
\begin{document}
\baselineskip=14pt minus 1pt

\parindent0pt
\maketitle
\begin{abstract}
\noindent 
We study unconditional subsequences of the canonical basis
$(\e_{rc})$ of elementary matrices in the Schatten class
$S^p$. They form the matrix counterpart to Rudin's $\Lambda(p)$
sets of integers in Fourier analysis.  In the case of $p$ an even integer, we find a sufficient
condition in terms of trails on a bipartite graph. We also
establish an optimal density condition and present a random
construction of bipartite graphs. As a byproduct, we get a new
proof for a theorem of Erd\H os on circuits in graphs.
\end{abstract}
\section{Introduction}

We study the following question on the Schatten class
$S^p$.\vskip\abovedisplayskip

\block{$(\dagger)$}{How many matrix coefficients of an
operator $x\in S^p$ must vanish so that the norm of $x$ has a
bounded variation if we change the sign of the remaining nonzero
matrix coefficients~?}\vskip\belowdisplayskip

Let $\Col$ be the set of columns and $\Row$ be the set of rows for
coordinates in the matrix, in general two copies of $\N$. Let
$I\se \Row\times \Col$ be the set of matrix coordinates of the
remaining nonzero matrix coefficients of $x$. Property
$(\dagger)$ means that the subsequence $(\e_{rc})_{(r,c)\in I}$
of the canonical basis of elementary matrices is an unconditional
basic sequence in 
$S^p$: $I$ forms a $\sigma(p)$ set in the terminology of
\cite[\S4]{ha98}.

It is natural to wonder about the operator valued case, where the
matrix coefficients are themselves operators in $S^p$. As the
proof of our main result carries over to that case, we shall
state it in the more general terms of complete $\sigma(p)$ sets.

We show that for our purpose, a set of matrix entries $I\se
\Row\times \Col$ is best understood as a bipartite graph. Its two vertex
classes are 
$\Col$ and $\Row$, whose elements will
respectively be termed ``column vertices'' and ``row
vertices''. Its edges join only row vertices $r\in \Row$ with column
vertices $c\in \Col$, this occurring exactly if $(r,c)\in I$.

We obtain a generic condition for $\sigma(p)$ sets in the case of
even $p$ (Th.~\ref{sig:thm}) that generalizes \cite[Prop.\
6.5]{ha98}. These sets reveal in fact as a matrix counterpart to
Rudin's $\Lambda(p)$ sets and we are able to transfer Rudin's
proof of
\cite[Th.~4.5(b)]{ru60} to a non-commutative context: his number
$r_s(E,n)$ is replaced by the numbers of Def.~\ref{def2}$(b)$
and we count trails between given vertices instead of
representations of an integer.

We also establish an upper bound for the intersection of a
$\sigma(p)$ set with a finite product set $\Row'\times \Col'$ (Th.\
\ref{dc}): this is a matrix  
counterpart to Rudin's \cite[Th.~3.5]{ru60}. In terms of
bipartite graphs, this intersection is the subgraph induced by
the vertex subclasses $\Col'\se\Col$ and $\Row'\se\Row$. 

The bound of Th.~\ref{dc} provides
together with Th.~\ref{sig:thm} a generalization of a theorem by Erd\H os
\cite[p.~33]{er64} on graphs without circuits of a given even
length. In the last part of this article, we present a random
construction of maximal $\sigma(p)$ sets for even integers $p$.

\vskipa {\bf Terminology } $\Col$ is the set
of columns and $\Row$ is the set of rows, in general both indexed
by $\N$. The set $V$ of all vertices is their disjoint union $\Row\amalg \Col$. An
edge on $V$ is a pair $\{v,w\}\se V$. A graph on $V$ is given by
its set of edges $E$. A bipartite graph on $V$ with vertex
classes $\Col$ and $\Row$ has only edges $\{r,c\}$ such that
$c\in\Col$ and $r\in\Row$ and may therefore be described
alternatively by the set $I=\{(r,c)\in\Row\times\Col:\{r,c\}\in
E\}$. A trail of length $s$ in a graph is a sequence
$(v_0,\dots,v_s)$ of $s+1$ vertices such that
$\{v_0,v_1\},\dots,\{v_{s-1},v_s\}$ are pairwise distinct edges
of the graph. A trail is a path if its vertices are pairwise
distinct. A circuit of length $p$ in a graph is a sequence $(v_1,\dots,v_p)$
of $p$ vertices such that
$\{v_0,v_1\},\dots,\{v_{p-1},v_p\},\{v_p,v_1\}$ are pairwise
distinct edges of the graph. A circuit is a cycle if its vertices
are pairwise distinct.

\vskipa {\bf Notation } $\T=\{z\in\C:|z|=1\}$. Let $q=(r,c)\in \Row\times
\Col$. The transpose of $q$ is $q^*=(c,r)$. The entry\index{entry}
(elementary matrix) $\e_q=\e_{rc}$ is the operator on $\ell_2$
that maps the $c$th basis vector on the $r$th basis vector and
all other basis vectors on $0$.  The matrix coefficient at
coordinate $q$ of an operator $x$ on $\ell_2$ is
$x_q=\trace\e_q^*x$ and its matrix representation is $(x_q)_{q\in
\Row\times \Col}=\sum_{q\in \Row\times \Col}x_q\e_q$.
The Schatten class $S^p$, $1\le p<\infty$, is the space of those
compact operators $x$ on $\ell_2$ such that
$\|x\|_p^p=\trace|x|^p=\trace(x^*x)^{p/2}<\infty$.  For $I\se
\Row\times \Col$, the entry space\index{entry subspace of $S^p$}
$S^p_I$ is the space of those $x\in S^p$ whose matrix
representation is supported by $I$: $x_q=0$ if $q\notin
I$. $S^p_I$ is also the closed subspace of $S^p$ spanned by
$(\e_q)_{q\in I}$.
The $S^p$-valued Schatten class $S^p(S^p)$ is the space of
those operators $x$ from $\ell_2$ to $S^p$ such that
$\|x\|_p^p=\trace(\trace|x|^p)<\infty$, where the inner trace is
the $S^p$-valued analogue of the usual trace.  The $S^p$-valued
entry space $S^p_I(S^p)$ is the
closed subspace spanned by the $x_q\e_q$ with $x_q\in S^p$ and
$q\in I$: $x_q=\trace\e_q^*x$ is the operator coefficient of $x$
at matrix coordinate $q$. Thus,
for even integers $p$ and $x=(x_q)_{q\in I}=\sum_{q\in I}x_q\e_q$
with $x_q\in S^p$ and $I$ finite, 
$$
\|x\|_p^p=\sum_{q_1,\dots q_p\in I}\trace x_{q_1}^*x_{q_2}\dots
x_{q_{p-1}}^*x_{q_p}\trace \e_{q_1}^*\e_{q_2}\dots
\e_{q_{p-1}}^*e_{q_p} .
$$

A Schur multiplier\index{Schur multiplier} $T$ on $S^p_I$
associated to $(\mu_q)_{q\in I}\in\C^I$ is a bounded operator on
$S^p_I$ such that $T\e_q=\mu_q\e_q$ for $q\in I$.  $T$ is
furthermore completely bounded\index{completely bounded Schur
multiplier} (c.b.~for short) if $T$ is bounded as the operator on
$S^p_I(S^p)$ defined by $T(x_q\e_q)=\mu_qx_q\e_q$
for $x_q\in S^p$ and $q\in I$.

We shall stick to this harmonic analysis type notation; let us
nevertheless show how these objects are termed with tensor
products: $S^p(S^p)$ is also $S^p(\ell_2\otimes_2\ell_2)$ 
endowed with $\|x\|_p^p=\trace\otimes\trace|x|^p$; one should write
$x_q\otimes\e_q$ instead of $x_q\e_q$; here 
$x_q=\Id_{S^p}\otimes\trace((\Id_{\ell_2}\otimes\e_q^*)x)$; $T$ is
c.b.~if $\Id_{S^p}\otimes T$ is bounded on
$S^p(\ell_2\otimes_2\ell_2)$. \vskipb

{\bf Acknowledgment } 
The first-named and last-named authors undertook this research
at the \'Equipe d'Analyse de l'Universit\'e Paris
6. It is their pleasure to acknowledge its kind  hospitality.
\section{Definitions and main results}
We use the notion of unconditionality in order to define the
matrix analogue of Rudin's ``commutative'' $\Lambda(p)$
sets.
\begin{dfn} 
Let $X$ be a Banach space. The sequence $(y_n)\se X$ is an
unconditional basic sequence in $X$ if there is a constant $D$
such that $$\Bigl\|\sum\vartheta_nc_ny_n\Bigr\|_X\le D\Bigl\|\sum
c_ny_n\Bigr\|_X$$ for every real \(\vs complex\) choice of signs
$\vartheta_n\in\{-1,1\}$
\(\vs $\vartheta_n\in\T$\)  
and every finitely supported sequence of scalar coefficients
$(c_n)$. The optimal $D$ is the real \(\vs complex\)
unconditionality constant of $(y_n)$ in $X$.
\end{dfn}
Real and complex unconditionality are isomorphically equivalent:
the complex unconditionality constant is at most $\pi/2$ times
the real one. The notions of unconditionality and multipliers are
intimately connected: we have
\begin{prp}
Let $(y_n)\se X$ be an unconditional basic sequence in $X$ and
let $Y$ be the closed subspace of $X$ spanned by $(y_n)$. The real
\(\vs complex\) unconditionality constant of $(y_n)$ in $X$ is
exactly the least upper bound for the norms
$\|T\|_{\script{L}(Y)}$, where $T$ is the multiplication operator
defined by $Ty_n=\mu_ny_n$, and the $\mu_n$ range over all real
\(\vs complex\) numbers with $|\mu_n|\le1$.
%
\end{prp}
Let us encompass the notions proposed in Question $(\dagger)$.

\begin{dfn}Let $I\se \Row\times \Col $ and $p>2$. 

$(a)$ \cite[Def.~4.1]{ha98} $I$ is a $\sigma(p)$ set if
$(\e_q)_{q\in I}$ is an unconditional basic sequence in
$S^p$. This amounts to the uniform boundedness of the family of
all relative Schur multipliers by signs
\begin{equation}\label{sigma:schur}
T_\vartheta\colon S^p_I\to S^p_I\ ,\ x=(x_q)_{q\in I}\mapsto
T_\vartheta x=(\vartheta_qx_q)_{q\in I}\hbox{ with }
\vartheta_q\in\{-1,1\}.
\end{equation}
By \cite[Lemma 0.5]{ha98}, this means that there is a constant
$D$ such that for every finitely supported operator
$x=(x_q)_{q\in I}=\sum_{q\in I}x_q\e_q$ with $x_q\in\C$
\begin{equation}\label{def:1}
D^{-1}\|x\|_p\le|\!|\!|x|\!|\!|_p\le\|x\|_p,
\end{equation}
where the second inequality is a convexity inequality that is
always satisfied \(see
\cite[Th.~8.9]{si80}\) and
\begin{equation}\label{def:2}
|\!|\!|x|\!|\!|_p^p=
\sum_c\Bigl(\sum_r|x_{rc}|^2\Bigr)^{p/2}
\vee
\sum_r\Bigl(\sum_c|x_{rc}|^2\Bigr)^{p/2}.
\end{equation}

$(b)$ \cite[Def.~4.4]{ha98} $I$ is a complete $\sigma(p)$ set if
the family of all relative Schur multipliers by signs
\Ref{sigma:schur} is uniformly c.b. By \cite[Lemma 0.5]{ha98},
$I$ is completely $\sigma(p)$ if and only if there is a constant
$D$ such that for every finitely supported operator valued
operator $x=(x_q)_{q\in I}=\sum_{q\in I}x_q\e_q$ with $x_q\in
S^p$
\begin{equation}\label{knc}
D^{-1}\|x\|_p\le|\!|\!|x|\!|\!|_p\le\|x\|_p,
\end{equation}
where the second inequality is a convexity inequality that is
always satisfied and $$
|\!|\!|x|\!|\!|_p^p=
\sum_c
\Bigl\|\Bigl(\sum_rx_{rc}^*x_{rc}\Bigr)^{1/2}\Bigr\|_p^p\vee
\sum_r
\Bigl\|\Bigl(\sum_cx_{rc}x_{rc}^*\Bigr)^{1/2}\Bigr\|_p^p.
$$
\end{dfn}
The notion of a complete $\sigma(p)$ set is stronger than that of
a $\sigma(p)$ set: Inequality \Ref{def:1} amounts to Inequality
\Ref{knc} tested on operators of the type 
$x=\sum_{q\in I}x_q\e_q$ with each $x_q$ acting on the same
one-dimensional subspace of $\ell_2$. It is an important open problem to
decide whether the notions differ. An affirmative answer would
solve Pisier's conjecture about completely bounded Schur
multipliers \cite[p.~113]{pi96}.

Notorious examples of $1$-unconditional basic sequences in all
Schatten classes $S^p$ are single columns, single
rows, single diagonals and single anti-diagonals --- and more
generally ``column'' (\vs ``row'') sets $I$ such 
that for each $(r,c)\in I$, no other
element of $I$ is in the column $c$ (\vs row $r$).

We shall try to express these notions in terms of trails on
bipartite graphs.  We proceed as announced in the Introduction:
then each example above is a union of disjoint star graphs in
which one vertex of one class is connected to some vertices of
the other class: trails in a star graph have at most length 2.
%
\begin{dfn}\label{def2}
Let $I\se \Row\times \Col $ and $s\ge1$ an integer. We consider $I$ as
a bipartite graph: its vertex set is $V=\Row\amalg \Col$ and its edge
set is $E=\bigl\{\{r,c\}\se V:(r,c)\in I\bigr\}$.

$(a)$ The sets of trails of length $s$ on the graph $I$ from the
column \(\vs row\) vertex $v_0$ to the vertex $v_s$ are
respectively
\begin{eqnarray*}
\script{\Col}^s(I;v_0,v_s)&=&\{(v_0,\dots,v_s)\in V^{s+1}:
v_0\in \Col\ \&\ \{v_i,v_{i+1}\}\in E\hbox{ pairwise distinct}\},\\
\script{\Row}^s(I;v_0,v_s)&=&\{(v_0,\dots,v_s)\in V^{s+1}:
v_0\in \Row\ \&\ \{v_i,v_{i+1}\}\in E\hbox{ pairwise distinct}\}.
\end{eqnarray*}
$(b)$ We define the Rudin numbers of trails starting respectively
with a column vertex and a row vertex by 
$c_s(I;v_0,v_s)=\#\script{\Col}^s(I;v_0,v_s)$ and 
$r_s(I;v_0,v_s)=\#\script{\Row}^s(I;v_0,v_s)$.
\end{dfn}
\rem
In other words, for an integer $l\ge1$,
\begin{eqnarray*}
c_{2l-1}(I;v_0,v_{2l-1})&=&\#\left[\vcenter{
\hbox{$(r_1,c_1),(r_1,c_2),(r_2,c_2),(r_2,c_3),\dots,(r_l,c_l)$}
\hbox{pairwise distinct in $I:c_1=v_0,r_l=v_{2l-1}$}}\right]\\
c_{2l}(I;v_0,v_{2l})&=&\#\left[\vcenter{
\hbox{$(r_1,c_1),(r_1,c_2),\dots,(r_l,c_l),(r_l,c_{l+1})$}
\hbox{pairwise distinct in $I:c_1=v_0,c_{l+1}=v_{2l}$}}\right]
\end{eqnarray*}
and similarly for $r_s(I;v_0,v_s)$.  If $s$ is odd, then
$c_s(I;v_0,v_s)=r_s(I;v_s,v_0)$ for all $(v_0,v_s)\in \Col\times
\Row$. But if $s$ is even, one Rudin number may be bounded while the
other is infinite: see
\cite[Rem.~6.4$(ii)$]{ha98}.\vskipb 
\section{$\sigma(p)$ sets as matrix $\Lambda(p)$ sets}
We claim the following result.
\begin{thm}\label{MAIN}
Let $I\se \Row\times \Col $ and $p=2s$ be an even integer.  If $I$ is a
union of sets $I_1,\dots,I_l$ such that one of the Rudin numbers
$c_s(I_j;v_0,v_s)$ or $r_s(I_j;v_0,v_s)$ is a bounded function of
$(v_0,v_s)$, for each $j$, then $I$ is a complete $\sigma(p)$
set.
\end{thm}
This follows from Theorem \ref{sig:thm} below: the union of two
complete $\sigma(p)$ sets is a complete $\sigma(p)$ set by
\cite[Rem.~after Def.~4.4]{ha98}; furthermore the transposed
set $I^*=\{q^*:q\in I\}\se \Col\times \Row$ is a complete $\sigma(p)$
set provided $I$ is.
Note that the case of $\sigma(\infty)$ sets (see \cite[Rem.\
4.6$(iii)$]{ha98}) provides evidence that Theorem \ref{MAIN}
might be a
characterization of complete $\sigma(p)$ sets for even $p$.
\begin{thm}\label{sig:thm}
Let $I\se \Row\times \Col $ and $p=2s$ be an even integer. If the Rudin
number $c_s(I;v_0,v_s)$ is a bounded function of $(v_0,v_s)$,
then $I$ is a complete $\sigma(p)$ set.
\end{thm}
This is proved for $p=4$ in \cite[Prop.~6.5]{ha98}. We wish to
emphasize that the proof below follows the scheme of the proof of
\cite[Th.~1.13]{ha98}. In particular, we make crucial use of
Pisier's idea to express repetitions by dependent Rademacher
variables (\cite[Prop.\ 
1.14]{ha98}).\vskipb

\dem
Let $x=\sum_{q\in I}x_q\e_q$ with $x_q\in S^p$. We have the
following expression for $\|x\|_p$.  $$
\|x\|_p^p=\trace\otimes\trace (x^*x)^s=\|y\|_2^2\quad
\hbox{with}\quad y=\overbrace{x^*xx^*\cdots x\pstar}^{s{\rm\ terms}},
$$
\ie $y$ is the product of $s$ terms which are alternatively $x^*$ and
$x$, and we set $x\pstar=x$ for even $s$, $x\pstar=x^*$ for odd
$s$.
Set $\Col\pstar=\Col$ for even $s$ and $\Col\pstar=\Row$ for odd $s$. Let
$(v_0,v_s)\in \Col\times \Col\pstar$ and
$y_{v_0v_s}=\trace\e_{v_0v_s}^*y$
be the matrix coefficient of $y$ at coordinate $(v_0,v_s)$. Then
we obtain by the rule of matrix multiplication
\begin{eqnarray*}
y=\sum_{q_1,\dots,q_s\in I}
(x_{q_1}^*\e_{q_1}^*)(x_{q_2}\e_{q_2})\dots
(x_{q_s}\pstar\e_{q_s}\pstar)
\end{eqnarray*}
\begin{eqnarray}\label{yvw}
y_{v_0v_s}=\sum_{
\scriptstyle (v_1,v_0),\,(v_1,v_2),\dots\,\in I
}
x_{v_1v_0}^*x_{v_1v_2}x_{v_3v_2}^*\dots x_{( v_{s-1},v_s
)\ppstar}\pstar.
\end{eqnarray}
Let $\script{E}$ be the set of equivalence relations on
$\{1,\dots,s\}$. Then
\begin{equation}\label{sim}
y=\sum_{\sim\in\script{E}}
\sum_{i\sim j\ssi q_i=q_j}
(x_{q_1}^*\e_{q_1}^*)(x_{q_2}\e_{q_2})\dots
(x_{q_s}\pstar\e_{q_s}\pstar).
\end{equation}
We shall bound the sum above in two steps.

$(a)$ Let $\sim$ be equality and consider the corresponding term
in the sum \Ref{sim}. The number of terms in the sum \Ref{yvw}
such that $\{v_{i-1},v_i\}\ne\{v_{j-1},v_j\}$ if $i\ne j$ is
$c_s(I;v_0,v_s)$. If $c$ is an upper bound for $c_s(I;v_0,v_s)$,
we have by the expression of the Hilbert--Schmidt norm and the
Arithmetic-Quadratic Mean Inequality
\begin{eqnarray*}
\lefteqn{\Bigl\|
\sum_{\scriptstyle q_1,\dots,q_s\atop\hbox{\scriptsize pairwise distinct}}
(x_{q_1}^*\e_{q_1}^*)(x_{q_2}\e_{q_2})\dots
(x_{q_s}\pstar\e_{q_s}\pstar)
\Bigr\|_2^2}\qquad\qquad&&\\
&=&\sum_ {( v_0,v_s)\in \Col\times \Col\ppstar}
\Bigl\|
\sum_{\scriptstyle v\in\script{\Col}^s(I;v_0,v_s)
}x_{v_1v_0}^*x_{v_1v_2}x_{v_3v_2}^*\dots
x_{(v_{s-1},v_s)\ppstar}\pstar
\Bigr\|_2^2
\\
&\le&c\sum_ {( v_0,v_s)\in \Col\times \Col\ppstar}
\sum_{\scriptstyle v\in\script{\Col}^s(I;v_0,v_s)
}\Bigl\| x_{v_1v_0}^*x_{v_1v_2}x_{v_3v_2}^*\dots
x_{(v_{s-1},v_s)\ppstar}\pstar
\Bigr\|_2^2
\\
&=&c\sum_{\scriptstyle q_1,\dots,q_s\atop\hbox{\scriptsize
pairwise distinct}}
\bigl\|(x_{q_1}^*\e_{q_1}^*)(x_{q_2}\e_{q_2})\dots(x_{q_s}\pstar\e_{q_s}\pstar)
\bigr\|_2^2\\
&\le&c\sum_{q_1,\dots,q_s}
\bigl\|(x_{q_1}^*\e_{q_1}^*)(x_{q_2}\e_{q_2})\dots(x_{q_s}\pstar\e_{q_s}\pstar)
\bigr\|_2^2\\
&=&c\Bigl\|\sum_{q_1,\dots,q_s}
|(x_{q_1}^*\e_{q_1}^*)(x_{q_2}\e_{q_2})\dots
(x_{q_s}\pstar\e_{q_s}\pstar)|^2
\Bigr\|_1
\end{eqnarray*}
Now this last expression may be bounded accordingly to
\cite[Cor.~0.9]{ha98} by
\begin{equation}\label{0.5}
c\Bigl(\Bigl\|\sum(x_q^*\e_q^*)(x_q\e_q)\Bigr\|_p\vee \Bigl\|\sum
(x_q\e_q)(x_q^*\e_q^*)\Bigr\|_p\Bigr)^p =c|\!|\!|x|\!|\!|_p^p:
\end{equation}
see \cite[Lemma 0.5]{ha98} for the last equality.

$(b)$ Let $\sim$ be distinct from equality. The corresponding
term in the sum \Ref{sim} cannot be bounded directly. Consider
instead $$
\Psi(\sim)=
\Bigl\|
\sum_{i\sim j\imp q_i=q_j}
(x_{q_1}^*\e_{q_1}^*)(x_{q_2}\e_{q_2})\dots
(x_{q_s}\pstar\e_{q_s}\pstar)
\Bigr\|_2
=\Bigl\|
\sum_{i\sim j\imp q_i=q_j}
\prod_{i=1}^sf_i(q_i)\Bigr\|_2
$$ with $f_i(q)=x_q\e_q$ for even $i$ and $f_i(q)=x_q^*\e_q^*$
for odd $i$. We may now apply Pisier's Lemma \cite[Prop.\
1.14]{ha98}: let $0\le r\le s-2$ be the number of one element
equivalence classes modulo $\sim$; then
\begin{equation}\label{PisHar}
\Psi(\sim)\le \|x\|_p^r(B|\!|\!|x|\!|\!|_p)^{s-r},
\end{equation}
where $B$ is the constant arising in Lust-Piquard's
non-commutative Khinchin inequality. In order to finish the
proof, one does an induction on the number 
of atoms of the partition induced by $\sim$, along the lines of
step 2 of the proof of \cite[Th.~1.13]{ha98}.\eck\vskip\belowdisplayskip

\rem
The Moebius inversion formula for partitions enabled Pisier
\cite{pi99} to
obtain the following explicit bounds in the computation above:
$$
\|y\|_2\le c^{1/2}|\!|\!|x|\!|\!|_p^s+
\sum_{0\le r\le s-2}{s\choose r}(s-r)!
\|x\|_p^r\bigl((3\pi/4)|\!|\!|x|\!|\!|_p\bigr)^{s-r}
$$ 
\begin{equation}\label{p-ortho}
\|x\|_p\le\bigl((4c)^{1/p}\vee 9\pi p/8\bigr)|\!|\!|x|\!|\!|_p.
\end{equation}
Let us also record the following consequence of his study of
$p$-orthogonal sums. The family $(x_q\e_q)_{q\in I}$
is $p$-orthogonal in the sense of \cite{pi99} if and only if the
graph associated to $I$ does not contain any circuit of length
$p$, so that we have by \cite[Th.\ 3.1]{pi99}:
\begin{thm}\label{piscycle}
Let $p\ge4$ be an even integer. If $I$ does not contain any circuit
of length $p$, then $I$ is a complete $\sigma(p)$ set with
constant at most $3\pi p/2$.
\end{thm}
%
\rem
Pisier proposed to us the following argument to deduce a weaker
version of Th.~\ref{sig:thm} from \cite[Th.~1.13]{ha98}. Let
$\Gamma=\T^V$ and $z_v$ denote the $v$th coordinate function on
$\Gamma$. Associate to $I$ the set $\Lambda=\{z_rz_c:(r,c)\in
I\}$. Let still $p=2s$ be an even integer. Then $I$ is a complete $\sigma(p)$ set if $\Lambda$ is a
complete $\Lambda(p)$ set as defined in \cite[Def.~1.5]{ha98},
which in turn holds if $\Lambda$ has property $Z(s)$ as given
in \cite[Def.~1.11]{ha98}. It turns out that this condition
implies the uniform boundedness of $$ c_t(I;v_0,v_t)\vee
r_t(I;v_0,v_t)\quad\hbox{for}\quad t\le s\ ,\ v_0,v_t\in V.  $$
For $p\ge8$, this implication is strict: in fact, the countable
union of disjoint cycles of length $4$ (``quadrilaterals'') $$
I=\bigcup\nolimits_{i\ge0}
\bigl\{(2i,2i),(2i,2i+1),(2i+1,2i+1),(2i+1,2i)\bigr\}
$$ satisfies $c_t(I;v_0,v_t)\vee r_t(I;v_0,v_t)\le2$ whereas
$\Lambda$ does not satisfy $Z(s)$ for any $s\ge4$.\vskipb

\rem\label{rem} 
This theorem is especially useful to construct c.b.~Schur
multipliers: by \cite[Rem.~4.6$(ii)$]{ha98}, if $I$ is a
complete $\sigma(p)$ set, there is a constant $D$ (the constant
$D$ in \Ref{knc}) such that for every sequence $(\mu_q)\in\C^{
\Row\times \Col }$ supported by $I$ and every operator
$T_\mu:(x_q)\mapsto(\mu_qx_q)$ we have $$
\|T_\mu\|_{\script{L}(S^p(S^p))}\le 
D\sup_{q\in I}|\mu_q|.  $$
\section{The intersection of a $\sigma(p)$ set with a finite
product set}
%
%
%
%
%
%
%
%
Let $I\se \Row\times \Col $ considered as a bipartite graph as in the
Introduction and let $I'\se I$ be the subgraph induced by the
vertex set $\Col'\amalg \Row'$, with $\Col'\se \Col$ a set of $m$ column
vertices and $\Row'\se \Row$ a set of $n$ row vertices. In other words,
$I'=I\cap \Row'\times \Col'$.  Let $d(v)$ be the degree of the vertex
$v\in \Col'\amalg \Row'$ in $I'$: in other words,
\begin{eqnarray*}
\forall c\in \Col'\quad d(c)&=&\#[I'\cap \Row'\times\{c\}],\\
\forall r\in \Row'\quad d(r)&=&\#[I'\cap\{r\}\times \Col'].
\end{eqnarray*}
Let us recall that the dual norm of \Ref{def:2} is $$
|\!|\!|x|\!|\!|_{p'}=\inf_
{\scriptstyle\alpha,\beta\in {S^{p'}}\atop\scriptstyle\alpha+\beta=x}
\biggl(\sum_c\Bigl(\sum_r|\alpha_{rc}|^2\Bigr)^{p'/2}\biggr)^{1/p'}+
\biggl(\sum_r\Bigl(\sum_c|\beta _{rc}|^2\Bigr)^{p'/2}\biggr)^{1/p'},
$$ where $p\ge2$ and $1/p+1/p'=1$ (see \cite[Rem.~after Lemma
0.5]{ha98}).
\begin{lem}\label{size:lem}
Let
$1\le p'\le2$ and $x=\sum_{q\in I'}x_q$. Then $$
|\!|\!|x|\!|\!|_{p'}^{p'}
\ge\doublesum_{(r,c)\in I'}\Bigl(\max\bigl(d(c),d(r)\bigr)^{1/2-1/p'}
|x_{rc}|\Bigr)^{p'} $$
\end{lem}
\dem
By the $p'$-Quadratic Mean Inequality and by Minkowski's
Inequality,
\begin{eqnarray*}
\lefteqn{
\biggl(\sum_{c\in \Col'}\Bigl(
\sum_{(r,c)\in I'}|\alpha_{rc}|^2\Bigr)^{p'/2}\biggr)^{1/p'}+
\biggl(\sum_{r\in \Row'}\Bigl(\sum_{(r,c)\in I'}|\beta_{rc}|^2\Bigr)^{p'/2}
\biggr)^{1/p'}}~&&\\
&\ge&
\Bigl(\sum_{c\in \Col'}d(c)^{p'/2-1}\sum_{(r,c)\in I'}|\alpha_{rc}|^{p'}
\Bigr)^{1/p'}+
\Bigl(\sum_{r\in \Row'}d(r)^{p'/2-1}\sum_{(r,c)\in I'}|\beta _{rc}|^{p'}
\Bigr)^{1/p'}\\
&\ge&
\Bigl(\doublesum_{(r,c)\in I'}
\bigl(d(c)^{1/2-1/p'}|\alpha_{rc}|+d(r)^{1/2-1/p'}|\beta_{rc}|\bigr)^{p'}
\Bigr)^{1/p'}
\end{eqnarray*}
The lemma follows by taking the infimum over all $\alpha,\beta$
with $\alpha_q+\beta_q=x_q$ for $q\in I'$ as one can suppose that
$\alpha_q=\beta_q=0$ if $q\notin I$; note further that
$1/2-1/p'\le0$.\eck
\begin{thm}\label{dc}
If $I$ is a $\sigma(p)$ set with constant $D$ as in \Ref{def:1},
then the size $\#I'$ of any subgraph $I'$ induced by $m$ column
vertices and $n$ row vertices, in other words the cardinal of any
subset $I'=I\cap \Row'\times \Col'$ with $\#\Col'=m$ and $\#\Row'=n$,
satisfies
\begin{eqnarray}
\#I'&\le&D^2\bigl(m^{1/p}n^{1/2}+m^{1/2}n^{1/p}\bigr)^2\label{dens}\\
\nonumber&\le&4D^2\min(m,n)^{2/p}\max(m,n).\end{eqnarray}
The exponents in this inequality are optimal even for a complete
$\sigma(p)$ set $I$ in the following cases:

$(a)$ if $m$ or $n$ is fixed;

$(b)$ if $p$ is an even integer and $m=n$.
\end{thm}
Bound \Ref{dens} holds \af if $I$ is a complete $\sigma(p)$
set. Density conditions thus do not so far permit to distinguish
$\sigma(p)$ sets and complete $\sigma(p)$ sets. One may
conjecture that Inequality \Ref{dens} is also optimal for $p$ not
an even integer and $m=n$: this would be a matrix counterpart to
Bourgain's theorem \cite{bo89} on maximal $\Lambda(p)$ sets. 
\vskipb

\dem
If \Ref{def:1} holds, then 
$\|x|_{I'}\|_p\le D|\!|\!|x|\!|\!|_{p}$ for all $x\in S^p$ by
Remark
\ref{rem} applied to $(\mu_q)$ the indicator function of $I'$, and
by duality $|\!|\!|x|_{I'}|\!|\!|_{p'}\le D\|x\|_{p'}$ for all
$x\in S^{p'}$ (compare with \cite[Rem.~4.6$(iv)$]{ha98}). Let
$$ y=\doublesum_{(r,c)\in I'}d(c)^{1/p'-1/2}\e_{rc}, $$ $$
z=\doublesum_{(r,c)\in I'}d(r)^{1/p'-1/2}\e_{rc}, $$ Then the $n$ rows
of $y$ are all equal, as well as the $m$ columns of $z$: $y$ and
$z$ have rank 1 
and a single singular value. By the norm inequality followed by
the $(2/p'-1)$-Arithmetic Mean Inequality,
\begin{eqnarray*}
\|y+z\|_{p'}&\le&\|y\|_{p'}+\|z\|_{p'}\\
&=&n^{1/2}\Bigl(\sum_{c\in \Col'}d(c)^{2/p'-1}\Bigr)^{1/2}+
m^{1/2}\Bigl(\sum_{r\in \Row'}d(r)^{2/p'-1}\Bigr)^{1/2}\\
&\le&n^{1/2}m^{1-1/p'}(\#I')^{1/p'-1/2}+
m^{1/2}n^{1-1/p'}(\#I')^{1/p'-1/2}.
\end{eqnarray*}
We used that $\sum_{c\in \Col'}d(c)=\sum_{r\in \Row'}d(r)=\#I'$. By
Lemma
\ref{size:lem} applied to $x=y+z$,
$$ (\#I')^{1/p'}\le
D(n^{1/2}m^{1-1/p'}+m^{1/2}n^{1-1/p'})(\#I')^{1/p'-1/2}, $$ and
we get therefore the first part of the theorem.

Let us show optimality in the given cases.

$(a)$ Suppose that $n$ is fixed and $\Col'=\Col$: $I'=\Row'\times \Col$ is a complete
$\sigma(p)$ set for any $p$ as a union of $n$ rows and
$\#I'=n\cdot m$.

$(b)$ is proved in \cite[Th.~4.8]{ha98}.\eck
\vskip\belowdisplayskip
\rem
If $n\nsim m$, the method used in \cite[Th.~4.8]{ha98} does not
provide optimal $\sigma(p)$ sets but the following lower bound.
Let $p=2s$ with $s\ge2$ an integer. Consider a prime $q$ and let
$k=s^{s-1}q^s$. By \cite[4.7]{ru60} and \cite[Th.~2.5]{ha98},
there is a subset $F\se\{0,\dots,k-1\}$ with $q$ elements whose
complete $\Lambda(2s)$ constant is independent of $q$. Let $m\ge
k$ and $0\le n\le m$ and consider the Hankel set $$
I=\bigl\{(r,c)\in\{0,\dots,n-1\}\times\{0,\dots,m-1\}:r+c\in
F+m-k\bigr\}.  $$ Then the complete $\sigma(p)$ constant of $I$
is independent of $q$ by \cite[Prop.~4.7]{ha98} and
$$\#I\ge\cases{ nq&if $n\le m-k+1$\cr
(m-k+1)q&if $n\ge m-k+1$.\cr} $$ If we choose $m=(s+1)k-1$, this
yields $$\#I\ge
{s^{1/s}\over(s+1)^{1+1/s}}\min(n,m)\max(m,n)^{1/s}.  $$
Random construction \ref{rc} provides bigger sets than this
deterministic construction; however, it also does not provide
sets that would show the optimality of Inequality \Ref{dens}
unless $s=2$.\vskipb

\section{Circuits in graphs}

Non-commutative methods yield a new proof to a theorem of Erd\H
os \cite[p.~33]{er64}. Note that its generalization by Bondy and
Simonovits \cite{bs74} is stronger than Th.~\ref{eee} below as it
deals with cycles instead of circuits. By Th.~\ref{piscycle} and
\Ref{dens}
\begin{thm}\label{eee}
Let $p\ge 4$ be an even integer. If $G$ is a nonempty graph
on $v$ vertices with $e$ edges without circuit of length $p$, then 
$$e\le18\pi^2 p^2\,v^{1+2/p}.$$
If $G$ is furthermore a bipartite graph whose two vertex classes
have respectively $m$ and $n$ elements, then
\begin{eqnarray}\label{luc}
e\le9\pi^2p^2\min(m,n)^{2/p}\max(m,n).
\end{eqnarray}
\end{thm}
\dem
For the first assertion, recall that a graph $G$ with $e$ edges
contains a bipartite subgraph with more than $e/2$ edges (see
\cite[p.~xvii]{bo78}). \eck\vskipb

\rem
\L uczak showed to us that \Ref{luc} cannot be optimal
if $m$ and $n$ are of very different order of magnitude. In
particular, let $p$ be a multiple of $4$. Let $e'$
be the maximal number of edges of a graph on $n$ vertices without
circuit of length $p/2$. If $m>pe'$, he shows that \Ref{luc} may
be replaced by $e<3m$.\vskipb
 
We also get the following result, which enables us to conjecture
a generalization of the theorems of Erd\H{o}s and Bondy and
Simonovits.
\begin{thm}\label{ggg}
Let $G$ be a nonempty graph on $v$ vertices with $e$ edges. Let
$s\ge 2$ be an integer.

$(i)$ If $$e>
8D^2\,v^{1+1/s}\quad\hbox{with }D>9\pi s/4,$$ then
one may choose two
vertices $v_0$ and $v_s$ such that $G$ contains more than
$D^{2s}/4$ pairwise distinct trails from $v_0$ to $v_s$, each of
length $s$ and with pairwise distinct edges.

$(ii)$ One may draw the same conclusion if $G$ is a bipartite
graph whose two vertex classes have 
respectively $m$ and $n$ elements and
$$e>4D^2\min(m,n)^{1/s}\max(m,n)\quad\hbox{with }D>9\pi s/4.$$ 
\end{thm}
\dem
$(i)$ According to \cite[p.~xvii]{bo78}, the graph $G$ contains
a bipartite subgraph with more than $e/2$ edges,
so that we may apply $(ii)$.

$(ii)$ Combining inequalities
\Ref{p-ortho} and \Ref{dens}, if $D>9\pi s/4$, then there are
vertices $v_0$ and $v_s$ such that the number $c$ of pairwise
distinct trails from $v_0$ to $v_s$, each of length $s$ and with
pairwise distinct edges, satisfies $(4c)^{1/2s}>D$.\eck
\vskipb

Two paths with equal endvertices are called independent if they
have only their endvertices in common.\vskipb

\quest Let $G$ be a graph on $v$ vertices with $e$ edges. Let
$s,l\ge 2$ be integers. Is it so that there is a constant $D$
such that if $e>Dv^{1+1/s}$, then $G$ contains $l$ pairwise
independent paths of length $s$ with equal endvertices~?\vskipb

\rem
Note that by Th.~\ref{dc}, the exponent $1+1/s$ is optimal in
Th.~\ref{ggg}$(i)$, 
whereas optimality of the exponent $1+2/p$ in
Th.~\ref{eee} is an 
important open question in Graph Theory (see \cite{luw99}).\vskipb

One may also formulate Th.~\ref{ggg}$(ii)$ in the following way.
\begin{thm}\label{fff}
%
If a bipartite graph $G_2(n,m)$ with $n$ and $m$ vertices in its
two classes avoids
any union of $c$ pairwise distinct trails along $s$ pairwise
distinct edges between two given vertices as a subgraph, where
the class of the first vertex is fixed, then the size $e$ of the
graph satisfies   
$$
e\le 4\max((4c)^{1/2s},9\pi s/4))\min(m,n)^{1/s}\max(m,n).
$$
\end{thm}
\section{A random construction of graphs}
Let us precise our construction of a random graph.
\begin{construction}\label{rc}
Let $\Col,\Row$ be two sets such that $\#\Col=m$ and $\#\Row=n$. Let
$0\le\alpha\le1$.  A random bipartite graph on $V=\Col\amalg \Row$ is
defined by selecting independently each edge in
$E=\bigl\{\{r,c\}\se V:(r,c)\in \Row\times \Col\bigr\}$ with the same
probability $\alpha$. The resulting random edge set is denoted by
$E'\se E$ and $I'\se \Row\times \Col$ denotes the associated random
subset.
\end{construction}
Our aim is to construct large sets while keeping down the Rudin
number $c_s$.
\begin{thm}
For each $\eps>0$ and for each integer $s\ge2$, there is an $\alpha$ such that Random
construction \ref{rc} yields subsets
$I'\se \Row\times \Col$ with size
$$\#I'\sim\min(m,n)^{1/2+1/s}\max(m,n)^{1/2-\eps}$$
and with $\sigma(2s)$ constant independent of $m$ and $n$ for $mn\to\infty$.
\end{thm}
\dem
Let us suppose without loss of generality that $m\ge n$. We want
to estimate the Rudin number of trails in $I'$.
Set ${\Col}\pstar=\Col$ for even $s$, ${\Col}\pstar=\Row$ for odd $s$ and let
$(v_0,v_s)\in \Col\times{\Col}\pstar$.
Let $l\ge1$ be a fixed integer. Then
\begin{eqnarray*}
\lefteqn{
\P{c_s(I';v_0,v_s)\ge l}=
\P{\exists\ l\hbox{ pairwise distinct trails }(v_0^j,\dots,v_s^j)\in
\script{\Col}^s(I';v_0,v_s)}}\qquad\qquad\\
&=&\P{E'\es\bigl\{\{v_{i-1}^j,v_i^j\}\bigr\}_{i,j}:
\{(v_0^j,\dots,v_s^j)\}_{j=1}^l\se\script{\Col}^s(\Row\times \Col;v_0,v_s)}\\
&\le&\sum_{k=\lceil l^{1/s}\rceil}^{ls}\#A_k\cdot\alpha^k,
\end{eqnarray*}
where $A_k$ is the following set of $l$-element subsets of trails
in $\script{\Col}^s(\Row\times \Col;v_0,v_s)$ built with $k$ pairwise
distinct edges $$
A_k=\Bigl\{\{(v_0^j,\dots,v_s^j)\}_{j=1}^l\se\script{\Col}^s(\Row\times
\Col;v_0,v_s):
\#\bigl\{\{v_{i-1}^j,v_i^j\}\bigr\}_{i,j}=k\Bigr\};
$$ the lower limit of summation is $\lceil l^{1/s}\rceil$
because one can build at most $k^s$ pairwise distinct trails of
length $s$ with $k$ pairwise distinct edges.

In order to estimate $\#A_k$, we now have to bound the number of
pairwise distinct vertices and the number of pairwise distinct
column vertices in each  set of $l$ trails $\{(v_0^j,\dots,v_s^j)\}_{j=1}^l\in A_k$. We claim that
\begin{eqnarray}\label{v1}
\#\{v_i^j:1\le i\le s-1,1\le j\le l\}
&\le&
k(s-1)/s,\\
\#\{v_{2i}^j:1\le i\le\lceil s/2 \rceil -1, 1\le j\le l\}
&\le& k/2.\label{v2}
\end{eqnarray}
%
The second estimate is trivial, because each column vertex
$v_{2i}^j$ accounts for two distinct edges
$\{v_{2i-1}^j,v_{2i}^j\}$ and $\{v_{2i}^j,v_{2i+1}^j\}$. For the
first estimate, note that each maximal sequence of $h$
consecutive pairwise distinct vertices
$(v_{a+1}^j,\dots,v_{a+h}^j)$ accounts for $h+1$ 
pairwise distinct edges $$\{v_a^j,v_{a+1}^j\}\ ,\ 
\{v_{a+1}^j,v_{a+2}^j\}\ ,\ \dots\ ,\ \{v_{a+h}^j,v_{a+h+1}^j\}\ ;$$
as $h\le s-1$, $h+1\ge hs/(s-1)$.
By \Ref{v1} and \Ref{v2}, $$ \#A_k\le
m^{k/2}n^{k/2-k/s}(k-k/s)^{ls-l}\le (ls)^{ls}m^{k/2}n^{k/2-k/s}:
$$ 
each element of $A_k$ is obtained by a choice of at most $k/2$
column vertices, a choice of at most $k/2-k/s$ row vertices, and
the choice of an arrangement with repetition of $ls-l$ out of at
most $k-k/s$ vertices.

Put
$\alpha=m^{-1/2}n^{-1/2+1/s}(\#\Col\cdot\#\Col\pstar)^{-\eps}$. Then
\begin{eqnarray*}
\P{\sup_{(v_0,v_s)
}c_s(I';v_0,v_s)\ge l}&\le&
\#\Col\cdot\#\Col\pstar\cdot(ls)^{ls}\sum_{k=\lceil
l^{1/s}\rceil}^{ls}(\#\Col\cdot\#\Col\pstar)^{-k\eps}\\
&\le& (ls)^{ls}
{(\#\Col\cdot\#\Col\pstar)^{1-\lceil
l^{1/s}\rceil\eps}\over1-(\#\Col\cdot\#\Col\pstar)^{-\eps}}.
\end{eqnarray*}
Choose $l$ such that $\lceil l^{1/s}\rceil\eps>1$. Then this
probability is little for $mn$ large. On the other hand, $\#I'$ is 
of order $mn\alpha$ with probability close to $1$.\eck\vskipb

\rem
This construction yields much better results for $s=2$. Keeping
the notation of the proof above and $m\ge n$, we get $k=2l$, $A_k={n\choose
l}$ and $$
\P{\sup_{(v_0,v_2)\in \Col\times \Col}c_2(I';v_0,v_2)\ge l}\le
m^2{n \choose l}\alpha^{2l}.$$ Let $l\ge2$ and
$\alpha=m^{-1/l}n^{-1/2}$. This yields sets $I'\se \Row\times \Col$
with size $$
\#I'\sim n^{1/2}m^{1-1/l}
$$ and with $\sigma(4)$ constant independent of $m$ and $n$.
This case has been extensively studied in Graph theory as the
``Zarankiewicz problem'': if $c_2(I';v_0,v_2)\le l$ for all
$v_0,v_2\in\Col$,
then the graph $I'$ does not contain a complete bipartite
subgraph on any two column vertices $v_0,v_2$ and $l+1$ row
vertices. Reiman (see \cite[Th.~VI.2.6]{bo78})
showed that then 
$$\#I'\le\bigl(lnm(m-1)+n^2/4\bigr)^{1/2}+n/2
\sim
l^{1/2}n^{1/2}m.$$
With use of finite projective geometries, he also showed that this
bound is optimal for 
$$n=l{q^{r+1}-1\over q^2-1}{q^r-1\over q-1} \quad,\quad
m={q^{r+1}-1\over q-1}$$
with $q$ a prime power and $r\ge2$ an
integer, and thus with $m\le n$: there seems to be no constructive
example of extremal graphs with $c_2(I';v_0,v_2)\le l$ and $m>n$
besides the trivial case of complete bipartite graphs with $m>n=l-1$.
\vskipb \rem In the case $s=3$, our result cannot be improved just by
refining the estimation of $\#A_k$. If we consider first $l$ distinct
paths that have their second vertex in common and then $l$ independent
paths, we get
$$\#A_{2l+1}\ge{m\choose l}n\quad,\quad \#A_{3l}\ge{m\choose
  l}{n\choose l}.$$
Therefore any choice of $\alpha$ as a monomial
$m^{-t}n^{-u}$ in the proof above must satisfy $t\ge(l+1)/(2l+1)$, $t+u\ge(2l+2)/(3l)$
and this yields sets with  
$$\#I'\preccurlyeq m^{1/2-1/2(4l+2)}n^{5/6-(7l+6)/(12l^2+6l)}.$$
%
%
%

\vskip\baselineskip
{\scriptsize
\halign{
#\hfill\quad&#\hfill\quad&#\hfill\quad\cr
Asma Harcharras&Stefan Neuwirth&Krzysztof Oleszkiewicz\cr
Mathematics Department&Laboratoire de Math\'ematiques&Institute of Mathematics\cr
University of Missouri-Columbia&Universit\'e de Franche-Comt\'e&Warsaw University\cr
Co\-lum\-bia MO 65211&25030~Besan\c con cedex&Banacha 2, 02-097 Warszawa\cr
U.S.A.&France&Poland\cr
harchars@math.missouri.edu&neuwirth@math.univ-fcomte.fr&koles@mimuw.edu.pl\cr
}}

\end{document}